%% file: biometrika.tex
\let\latexarabic\arabic
\let\latexdocument\document
\let\latexenddocument\enddocument

\RequirePackage[thmmarks]{ntheorem}

\makeatletter
\renewtheoremstyle{plain} 
  {\item[\hskip\labelsep \theorem@headerfont ##1\ \textup{##2}\theorem@separator]} 
  {\item[\hskip\labelsep \theorem@headerfont ##1\ \textup{##2}\ (##3)\theorem@separator]}
\makeatother

\documentclass[article,lineno]{biometrika}

\let\document\latexdocument
\let\enddocument\latexenddocument
\AtEndDocument{\printhistory}
\let\arabic\latexarabic
\def\rm{}

\usepackage{amsfonts, amsmath, amssymb, bbm}

\usepackage{times}
\usepackage{bm}
\usepackage{moreverb,url}
\usepackage{multirow}
\usepackage{float}
\usepackage{tikz}
\usetikzlibrary{arrows,shapes.arrows,shapes.geometric,shapes.multipart, decorations.pathmorphing,positioning,shapes.swigs}

\usepackage{mathtools}
\mathtoolsset{showonlyrefs}
\graphicspath{{./Figures/}}
\usepackage{color}
\definecolor{darkred}{RGB}{100,0,0}
\definecolor{darkgreen}{RGB}{0,100,0}
\definecolor{darkblue}{RGB}{0,0,150}

\definecolor{purple}{rgb}{0.4,.1,.9}

\usepackage{hyperref}
\hypersetup{colorlinks=true, linkcolor=darkred, citecolor=darkgreen, urlcolor=darkblue}

\usepackage[plain,noend]{algorithm2e}

\makeatletter
\renewcommand{\algocf@captiontext}[2]{\quad #1\algocf@typo. \AlCapFnt{}#2} 
\def\@algocf@capt@plain{top}
\renewcommand{\algocf@makecaption}[2]{%
  \addtolength{\hsize}{\algomargin}%
  \sbox\@tempboxa{\algocf@captiontext{#1}{#2}}%
  \ifdim\wd\@tempboxa >\hsize
    \hskip .5\algomargin%
    \parbox[t]{\hsize}{\algocf@captiontext{#1}{#2}}
  \else%
    \global\@minipagefalse%
    \hbox to\hsize{\box\@tempboxa}
  \fi%
  \addtolength{\hsize}{-\algomargin}%
}
\makeatother

\renewcommand{\P}{\operatorname{\mathbb{P}}}

\begin{document}

\jname{Biometrika}
\jyear{2021}
\jvol{103}
\jnum{1}
\accessdate{Advance Access publication on 31 August 2022}

\received{{\rm 2} January {\rm 2017}}
\revised{{\rm 8} June {\rm 2021}}


\markboth{Andrew Ying, Ronghui Xu}{Miscellanea}

\title{On Defense of the Hazard Ratio}

\author{Andrew Ying}
\affil{\email{aying9339@gmail.com}}

\author{Ronghui Xu}
\affil{Herbert Wertheim School of Public Health, Department of Mathematics, and Halicioglu Data Science Institute, University of California San Diego,\\ La Jolla, California 92093, U.S.A.
	\email{rxu@ucsd.edu }}

\maketitle
\begin{abstract}
In this short communication, we describe the recent debate on whether the hazard function should be used for causal inference in time-to-event studies and consider three different potential outcomes frameworks (by Rubin, Robins, and Pearl, respectively) as well as use the single-world intervention graph to show mathematically that the hazard function has causal interpretations 
 under all three frameworks. 
In addition, we argue that the hazard ratio over time can provide a useful interpretation in practical settings. 
\end{abstract}
\begin{keywords}
Causal Inference; Cox Model; Hazard Function; Potential Outcomes; Single-World Intervention Graph (SWIG); Survival Analysis. 
\end{keywords} 

\section{Background}

Modeling the hazard ratios through the Cox model is perhaps one of the most celebrated and best-adopted approaches for analyzing time-to-event data, partly due to its flexible semiparametric modeling form  \citep{reid:94}. It has been commonly used to analyze randomized clinical trials with survival endpoints. 
On causal interpretation of the hazard function, \citet{hernan2010hazards} first ``blew the whistle.'' 
He warned that a time-fixed hazard ratio can be misleading and a time-varying one has no causal interpretation, mainly because there is selection bias caused by differential survival distributions between the groups unless the null hypothesis of no treatment effect holds. He claimed that the selection bias is present even in the absence of unmeasured confounding, measurement error, and model misspecification. He illustrated the selection bias via an example from the Women’s Health Initiative \citep[WHI]{anderson2004effects, prentice2005statistical}.
Following that, \citet{aalen2015does} formalized this selection bias as a  collider bias when conditioning on a risk set in the presence of heterogeneity.
More recently \citet{martinussen2020subtleties} attempted 
to provide more insight into the subtle interpretation of hazard contrasts,  
as well as constructed an alternative estimand called the ``causal hazard ratio'' in order to deliver a causal interpretation.

\citet{prentice2022intention} presented a defense of hazard rate 
modeling, and Cox regression in particular, mainly from the applications perspective, for the intent-to-treat reporting of causal effects in randomized controlled trials. They elaborated with a comprehensive survival analysis of the massive Women’s Health Initiative's randomized, placebo-controlled hormone replacement therapy trials. 

In this paper, we formalize a mathematical defense for the hazard ratios. We prove that there is indeed no selection bias and it is safe to interpret hazard ratios causally. We refute the previous claims point by point below. Finally, we corroborate \citet{prentice2022intention}'s defense through more applications.

For the rest of this section we introduce some notation. Define $T$ as the time to event of interest, $A$  a randomized binary treatment assignment, and $L$  baseline covariates that can be partially or completely unobserved. We assume no censoring here in order to keep the focus on the hazard function itself. 
Denote $Y(t) = \mathbbm{1}(T \geq t)$ the at-risk process at time $t$, indicating whether a subject has survived at least to time $t$. The hazard function for $T$ is defined as 
\begin{equation}
    \lambda(t|\cdot) := \lim_{\Delta t \to 0+}\frac{1}{\Delta t}\P(t \leq T < t + \Delta t|T \geq t, \cdot) = \lim_{\Delta t \to 0+}\frac{1}{\Delta t}\P( Y(t + \Delta t) = 0|Y(t) = 1, \cdot),
\end{equation}
where $\P$ denotes the underlying probability measure. 

Analyses of time-to-event endpoints in randomised experiments are commonly based on the Cox proportional hazards model
\begin{equation}\label{eq:fixed}
    \lambda(t|A) = \lambda_0(t)\exp(\beta A),
\end{equation}
for a binary treatment $A$. Define also the following saturated generalized Cox model allowing time-varying treatment effect:
\begin{equation}\label{eq:varying}
    \lambda(t|A) = \lambda_0(t)\exp\{\beta(t) A\},
\end{equation}
so that 
\begin{equation}
    \exp\{\beta(t)\} = \frac{\lambda(t|A = 1)}{\lambda(t|A = 0)}. 
\end{equation}
We will discuss whether and when $\beta$ and $\beta(t)$ can be causallly interpreted. To that end, we introduce the potential outcomes \citep{rubin1974estimating, holland1986statistics}. Denote $T_a$ the potential time to event if $A$ were to be set to $a$, and correspondingly $Y_a(t) = \mathbbm{1}(T_a \geq t)$. Since $A$ is randomized, we have
\begin{equation}\label{eq:rand}
    A \perp (T_1, T_0, Y_1(t), Y_0(t),  L).
\end{equation}
We assume the standard causal consistency assumption
\begin{equation}
    T = T_A,
\end{equation}
and the positivity assumption that 
\begin{equation}
    0 < \P(A = 1) < 1.
\end{equation}

\section{Literature review}\label{sec:review}

\subsection{\citet{hernan2010hazards} stated that a time-fixed hazard ratio is misleading
}

\citet{hernan2010hazards} claimed that although the hazard ratios
may change over time in studies, often in practice a single hazard ratio is averaged over the duration of the study’s follow-up. As a result, the conclusions from the study may critically depend on the duration of the follow-up.

\subsection{\citet{hernan2010hazards}, \citet{aalen2015does} and \citet{martinussen2020subtleties}
claimed that a time-varying hazard ratio has no causal interpretation}\label{sec:timevarying}
Meanwhile by \eqref{eq:rand} and consistency,
\begin{equation}\label{eq:counterhr}
    \exp\{\beta(t)\} = \lim_{dt \to 0+}\frac{\P(t \leq T < t + dt|T \geq t, A = 1)}{\P(t \leq T < t + dt|T \geq t, A = 0)} = \lim_{dt \to 0+}\frac{\P(Y_1(t + dt) = 0|Y_1(t) = 1)}{\P(Y_0(t + dt) = 0|Y_0(t) = 1) }.
\end{equation}
\citet{hernan2010hazards}, \citet{aalen2015does} and \citet{martinussen2020subtleties}  concluded that \eqref{eq:counterhr} cannot be causally interpreted because the risk sets, i.e.~the conditioning sets comprised of the subsets of individuals who have not previously failed, differ beyond the first event time. 

In addition, \citet{hernan2010hazards} stated that ``differential selection of less susceptible women over time ... is the built-in selection bias of period-specific hazard ratios.'' 



\subsection{\citet{aalen2015does} claimed that the selection bias results from a collider bias}

\begin{figure}[ht!]
    \centering
    \input{Figures/DAG}
    \caption{A directed acyclic graph describing the data generating process.}
    \label{fig:dag}
\end{figure}
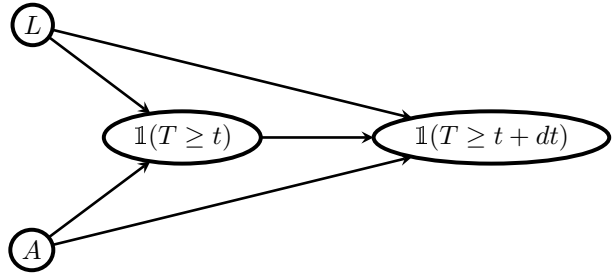

\citet{aalen2015does} provided a visualization to illustrate the selection bias, which we recreate here in Figure \ref{fig:dag}, as a causal directed acyclic graph (DAG) describing the data generating process of the observed data. 
In Figure \ref{fig:dag}, both $A$ and $L$ point to the risk sets $Y(t)$ and $Y(t + dt)$ admitting that $A$ and $L$ affect $Y(t)$ and $Y(t + dt)$ possibly; there is no arrow between $A$ and $L$ because $A$ is randomized. 
Since $A$ is randomized, a direct 
comparison of $Y(t)$ for any $t$ (including $t + dt$) between the treatments groups 
is a valid assessment of the causal effect. 
On the other hand, the node $Y(t) = \mathbbm{1}(T \geq t)$ is a collider and if one considers the probability of surviving up to time $t + dt$ conditional on survival up to time $t$, then 
the non-causal path $A \to Y(t) \leftarrow L \to Y(t + dt)$ is activated. 
This path cannot be closed unless $L$ is completely known, which implies that we generally have $A \not\perp L | Y(t) = 1$.  This can also been seen intuitively  if the treatment has any non-zero effect, since at any $t>0$ the two groups would have failed at different rates and are no longer exchangeable, despite the fact that  
the treatment was randomly assigned at $t = 0$. 

We note that the above DAG and argument are based on the observed variables instead of the potential variables. And the risk set, or equivalently, $Y(t)$, is also observed as opposed to being potential.

\subsection{\citet{martinussen2020subtleties} proposed an alternative ``causal'' HR}

Given their conclusion that the hazard ratio does not deliver causal interpretation,  \citet{martinussen2020subtleties} proposed the so-called ``causal hazard ratio, ''
\begin{equation}\label{eq:causalhr}
    \frac{\lim_{dt \to 0+}\P(t \leq T_1 < t + dt|T_1 \geq t, T_0 \geq t)/dt}{\lim_{dt \to 0+}\P(t \leq T_0 < t + dt|T_1 \geq t, T_0 \geq t)/dt}.
\end{equation}
\citet{martinussen2020subtleties} believed that \eqref{eq:causalhr} delivers valid causal interpretation because it compares quantities conditioned on the same subpopulation.

\section{Defense of hazard ratio in a causal setting}

\subsection{Time-fixed hazard ratio delivers causal interpretation}

We first confirm that a time-fixed hazard ratio $\beta$ delivers causal interpretation. Note that under model \eqref{eq:fixed}, consistency, and positivity,
\begin{equation}
    \exp(\beta) = \frac{\log \P(T_1 > t)}{\log \P(T_0 > t)}.
\end{equation}
Indeed this is confirmed in \citet{martinussen2020subtleties}.

While it might be the case that the underlying true hazard ratio changes over time, any practitioner understand the need for 
parsimony. 
\citet{prentice2022intention} showed that among other things, a single averaged hazard ratio is more sensitive to early treatment effect than, for example, the restricted mean survival time. 

\subsection{Time-varying hazard ratio has a causal interpretation}\label{sec:frameworks}

We note that the definitions of causal interpretation vary across different frameworks of causality. We examine below three different 
frameworks.
\begin{enumerate}
    \item[(a)] {\it Rubin's framework}:

The definition of causal interpretation in \citet{hernan2010hazards}, \citet{aalen2015does} and  \citet{martinussen2020subtleties} is best formalized from Donald Rubin's framework. According to definitions found in, for example, \cite{rubin1974estimating, rubin1978bayesian} and \citet{frangakis2002principal}, the above contrast \eqref{eq:counterhr} does not convey any causal interpretation. Given a random sample indexed by units $i$, according to \citet[Equation (2.1)]{frangakis2002principal}, in our notation: ``a causal effect of assignment on the outcome $T$ is defined to be a comparison between the potential outcomes on a common set of units, e.g., a comparison between $\{T_{1, i}:i \in \text{set}_1\}$ and $\{T_{0, i}: i \in \text{set}_0\}$ given the groups of units, $\text{set}_1$ and $\text{set}_0$, being compared are identical. '' 
Under this definition, $\beta(t)$ cannot be interpreted causally because the numerator and denominator in \eqref{eq:counterhr} are conditioned on $\text{set}_1 = \{i: Y_{1, i}(t) = 1\}$ and $\text{set}_0 = \{i: Y_{0, i}(t) = 1\}$, which cannot be identical except for under the sharp null that there is no individual treatment effect. 

Despite of the above, we can still rewrite
\begin{equation}
    \exp\{\beta(t)\} = \frac{\lim_{dt \to 0+}\log \P(T_1 > t + dt|T_1 > t)/dt}{\lim_{dt \to 0+}\log \P(T_0 > t + dt|T_0 > t)/dt}
    = \frac{ d \log \P(T_1 > t )/dt}{ d \log \P(T_0 > t )/dt},
\end{equation}
and in this way the time-varying hazard ratio seems to align with Rubin's framework.

\vskip .1in
    \item[(b)] {\it Robins' framework}:

\citet[Page 7, Technical Point 1.1]{hernan2020causal} defined: ``$\cdots$population causal effect may also be defined as a contrast of, say, medians, variances, hazards, or cdfs of counterfactual outcomes. In general, a population causal effect can be defined as a contrast of any functional of the marginal distributions of counterfactual outcomes under different actions or treatment values. ''
Since the hazard function itself is a functional of the marginal distributions, and hazard ratio is a contrast, the hazard ratio should be causally interpretable in this sense.

\vskip .1in
\item[(c)] {\it Pearl's framework}:

\citet[Page 70, Definition 3.2.1]{pearl2009causality} defined: ``Given two disjoint sets of variables, $X$ and $Y$, the causal effect of $X$ on $Y$, denoted either as $\P(y|\hat x)$ or as $\P(y|\text{do}(x))$, is a function from $X$ to the space of probability distributions on $Y$.''
The latter translates to potential outcomes with our notation as given 
$A$ and $T$, the causal effect of $A$ on $T$, denoted as $\P(T_a)$, is a function from $A$ to the space of probability distributions on $T$. Note that the knowledge of $\P(T_a)$ is equivalent to knowing the counterfactual hazard function $\lim_{dt \to 0+}\P\{Y_a(t + dt) = 0|Y_a(t) = 1\}/dt$ throughout time $t$. 
Therefore contrasting hazard functions between the treatment groups should yield valid causal interpretation under Judea Pearl's framework.
\end{enumerate}


Another confusion point in \citet{hernan2010hazards}  is about `less susceptible women'. This has been well understood in the literature as unobserved heterogeneity, i.e.~frailty, since at least as early as \cite{Lancaster}. It has also been shown in \cite{omor:john} that denoting $V$ the multiplicative frailty on the hazard function,  $E(V|T>t)$ is non-increasing in $t>0$, with no distributional assumption on $V$ required. In other words, the `differential selection ... over time' is due to the existence of unobserved heterogeneity, also referred to as over-dispersion, and not `period-specific hazard ratios' as speculated in \citet{hernan2010hazards}.

\subsection{No collider bias}

\begin{figure}[H]
    \centering
    \input{Figures/SWIG}
    \caption{A single-world intervention graph (SWIG) describing the data generating process for the potential outcomes.}
    \label{fig:swig}
\end{figure}
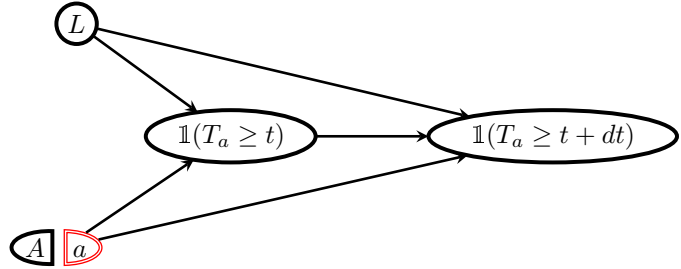

Here as opposed to \citet{aalen2015does} we draw a single-world intervention graph (SWIG) \citep{richardson2013single} in Figure \ref{fig:swig} which depicts the data generating process like the DAG in Figure \ref{fig:dag}, except on the potential outcomes. From a SWIG, one can better read dependency among variables in the counterfactual world by d-separation. To form a SWIG from a DAG, one first copies the DAG including all nodes and arrows. One then splits the treatment variable node $A$ into a random part $A$ and an intervention part $a$. All arrows pointing to $A$ in the original DAGs still point to $A$ in the SWIG. All arrows pointing to others from $A$ in the original DAGs now are pointing to others from the intervention node $a$ in the SWIG. Lastly, one turns all the post-treatment variables into potential variables. Following this way, the DAG in Figure \ref{fig:dag} is transformed into the SWIG in Figure \ref{fig:swig}. Since there is no arrow between $A$ and $Y_a(t)$, one can read that $A \perp L | Y_a(t) = 1$ and hence \emph{$Y_a(t)$ is not a collider}. It follows that the compositions of the groups of treated and untreated survivors at time $t$ are exchangeable in the counterfactual world, therefore the comparison in \eqref{eq:counterhr} still delivers a causal meaning.

The above shows that the hazard function is in fact conditional upon the counterfactual at-risk process $Y_a(t)$ under randomization, 
instead of the observed at-risk process $Y(t)$. Such an observation equips $\beta(t)$ with a causal interpretation because $A$ and $Y_a(t)$ are independent according to \eqref{eq:rand}.

\subsection{Causal hazard ratio is not satisfactory}

According to Section \ref{sec:frameworks}, the causal hazard ratio as defined in \eqref{eq:causalhr} can be causally interpreted by Rubin's causal framework because the numerator and the denominators condition on the same subpopulation. However, it does not fit Robins' framework because it is not a contrast of any functional of the marginal distributions of counterfactual outcomes, but rather the joint distributions. It also does not fit Pearl's framework for the same reason. Furthermore, it can never be nonparametrically identified because it hinges upon the joint distribution of $(T_1, T_0)$, which cannot be simultaneously observed. \cite{axel:nevo} developed an approach for sensitivity analysis, which unfortunately cannot be used as primary analysis of any randomized clinical trials. Last but not least, the causal hazard ratio is a type of ``survivor average causal effect," of which
\citet[Section 7]{dawid2000causal} questioned the use and called it ``a fundamental use of fatalism.'' He stated ``... it is only under the unrealistic assumption of fatalism that this group has any meaningful identity, and thus only in this case could such inferences even begin to have any useful content.'' 

In our own experience analyzing observational studies of exposure during  pregnancy, while live birth is a post treatment outcome and should be used to form survivor principal strata \citep{ying2020causal}, in practice it is nonetheless standard to stratify by the observed live birth outcome 
\citep{chambers2001postmarketing, cham:etal:2010, cham:2013, cham:2016, cham:etal:2019, cham:plusone, cham:hydroxy}. This is how the primary outcome results are presented in the teratology and dysmorphology literature, as well as the United States Food and Drug Administration vaccine and medication labelling.  

\section{hazard ratio in more complex settings}

As mentioned earlier, the Women's Health Initiative examples described in \citet{prentice2022intention} provided, in our view, convincing applications of hazard ratio in the past several decades in medical research. Since the publication of the WHI results, following the drastic change in clinical practice for postmenopausal women, there was soon thereafter report of reduction in postmenopausal breast cancer incidence in the US and elsewhere, and reductions in the US health care costs from the reduction in the use of the hormone therapies were estimated at \$37.1 one decade later. The early stopping of both of the WHI trials has to do with the fact that, without making the proportional hazards assumption, the average hazard ratio  is more sensitive to detecting early differences than the alternative restricted mean survival time, for example  \citep{prentice2022intention}. 

In the following, we provide other examples to illustrate the desired interpretation of $\beta(t)$, and more generally, changing hazard ratios which can be due to changing treatments in practice. 

In a two-group setting where the proportional hazards assumption is violated and the hazard functions cross each other, a representative example is organ transplantation \citep{oq:pess}. The picture that was shown in \cite{oq:pess} is in fact crossing survival curves, and this is when one wants to know, which treatment is better. In such a situation, the physician and their patient would  like to know, that compared to the same patient without a transplant, the patient who receives it is initially at a higher risk immediately following the transplant surgery but, once the patient survives the surgery and also does not develop the so-called graft-versus-host disease (GVHD), in the long run, the patient is much better off (i.e.~with lower hazard for death etc.) compared to the counterfactual scenario where they did not receive the transplant. This is different from the hypothetical situation where 
the procedure of organ transplantation has simply ``selected" those patients who are prone to developing GVHD and ``killed" them first, leaving the rest of the patients to live a much longer life but with no real benefit of getting a new organ. 
Both situations will appear as initially higher hazard then lower hazard compared to the counterfactual one, so
how can we tell if there is real treatment effect and not just the `selection' effect? 
One answer would be to look at the survival curves; in the latter situation the treated survival curve will remain lower than the untreated one, and they will not cross. 
{(One may also devise some analysis to separate the effects if possible, otherwise sensitivity analysis but in reality there has not been such a need, as science has obviously found a way to discover the benefit of transplant.)}
No matter which solution we adopt, however, in understanding and explaining what is going on, we still need to resort to the concept of hazard function and describe it as: compared to the same patient without a transplant, the patient who receives it is initially at a higher (instantaneous) risk immediately following the transplant surgery but, once the patient survives the surgery and also does not develop the GVHD, in the long run the patient is much better off \emph{due to the benefit of the new organ}.
The obvious common-sense lesson here is, unlike those who caution against the use of time-varying $\beta(t)$ at some later time point $t$ and, as a result mistakenly concluding that the treatment is beneficial in the second hypothetical situation above, 
one should at least always look at the whole curve of $\beta(t)$, and not just some later time point $t$ in isolation. 


The above leads to a very simple concept that is well understood in physics: distance versus velocity versus acceleration. Or equivalently, the concept of derivative in calculus, and its applications in science. Under smoothness assumptions, knowing one we know the other. This then leads naturally to dynamic systems where derivatives are useful, and is increasingly being applied to treatment processes. In the context of dynamic treatment regimes, it is conceivable that the treatment immediately impacts the hazard, and this impact is then reflected in cumulative quantities such as the survival probabilities. 
As a very simple example of time to exposure \citep[or treatment discontinuation]{yang2018modeling}, denote $a_d(t) = I(d\leq t)$ the treatment process if a subject becomes exposed at time $d>0$. 
A dynamic regime marginal structural Cox model in this case can be written as 
\begin{equation}
\lambda_{d}(t) = \lambda_\infty(t)\exp\{\beta a_d(t)\}, \label{eq:dr-MSM_Cox}
\end{equation}
where 
$\lambda_d$ is the hazard function of the potential survival outcome $T_d$, and $\lambda_\infty(t)$ is the hazard function of $T_\infty$, i.e.~for someone who is never exposed. 
Clearly, the hazard for someone who is not yet exposed is the same as $\lambda_\infty(t)$, until the subject becomes exposed, and then the hazard ratio becomes $ \exp(\beta) $.
It can be shown that for $t\leq d$, the survival probability 
$    S_{d}(t) = P(T_d \geq t) = S_\infty(t)$;
while for $t>d$, 
$    S_{d}(t)  
= S_\infty(d) \left\{ {S_\infty(t)}/{S_\infty(d)} \right\}^{\exp(\beta)}     $. 

As another example, during the recent COVID-19 pandemic, vaccine effectiveness is measured as one minus the hazard ratio.  \cite{Lin:etal:2022, lin2022effects} 
used a time-varying hazard ratio, by placing a change point every month, to estimate the \emph{current risk} of COVID-19 after receiving the vaccine. Because the hazard pertains to the current risk, estimation of the hazard ratio from the time of the first dose does not bias the estimation of the hazard ratio since full vaccination (changing treatment) and enables evaluation of the vaccine effect during the ramp-up period. We refer interested readers to the publications for more details. 

A final example is multi-state systems. The simplest multi-state system is a three-state system, that can be applied to semi-competing risks \citep{xu:kalb:2010, zhang:yiran}. A multi-state system is completely characterized by its transition intensities, also called the transition hazards \citep{geskus2015},  which are the same as the hazard function in a two-state, i.e.~time-to-event setting discussed in this paper. As stated above, any treatment effect would first modify the transition hazard(s), then subsequently any cumulative risks one might wish to consider. As shown in \cite{zhang:yiran}, expanding the concept of potential outcomes, we would consider \emph{potential systems} that consist of multiple states and transition times between them. 
In the simplest scenario of a binary point treatment at baseline, these systems correspond to two potential worlds, where the variables evolve over time. It is inevitable then, if there is a treatment effect, that the two potential worlds will be different after time zero, but that should not prevent us from looking at the `velocity' or even 'acceleration' by which they change along the time axis.

\section{Conclusion}
Previous literature critical of the hazard ratio has attempted to unduly constrain its causal interpretation mathematically. In this short communication, we have clarified that the hazard function has causal interpretations under all three potential outcomes frameworks commonly adopted for causal inference. In practice, the measure used to quantify treatment effects depends on the application at hand. The hazard ratio has been successful in the past in studying interventions for cancer, HIV AIDS, and many other diseases. We hope that a debate like this will lead to further understanding and sometimes, to rediscovering what we have already learned from the past. 





\section*{Acknowledgement}

The authors would like to acknowledge the Joint Statistical Meetings 2022 round table discussion sponsored by the Statistics in Epidemiology Section, which motivated the writing of this short communication. The authors also acknowledge discussion at the Lifetime Data Science (LiDS) conference 2023 which furthered our thinking on the subject.

\bibliographystyle{biometrika}
\bibliography{ref}

\end{document}

%% file: Figures/DAG.tex
\begin{tikzpicture}

\tikzset{line width=1.5pt, outer sep=0pt,
ell/.style={draw,fill=white, inner sep=2pt,
line width=1.5pt},
swig vsplit={gap=5pt, line color right=red,
inner line width right=0.5pt}};

\node[name=L, ell, shape=ellipse] {$L$};

\node[name=Yt, below right=15mm of L, ell, shape=ellipse]{$\mathbbm{1}(T \geq t)$};

\node[name=Ytd, right=15mm of Yt, ell, shape=ellipse]{$\mathbbm{1}(T \geq t + dt)$};

\node[name=A, below left =15mm of Yt, ell, shape=ellipse]{$A$};

\draw[->,line width=1.0pt,>=stealth](L) to (Yt);
\draw[->,line width=1.0pt,>=stealth](L) to (Ytd);
\draw[->,line width=1.0pt,>=stealth](Yt) to (Ytd);
\draw[->,line width=1.0pt,>=stealth](A) to (Yt);
\draw[->,line width=1.0pt,>=stealth](A) to (Ytd);

\end{tikzpicture}

%% file: Figures/SWIG.tex
\begin{tikzpicture}

\tikzset{line width=1.5pt, outer sep=0pt,
ell/.style={draw,fill=white, inner sep=2pt,
line width=1.5pt},
swig vsplit={gap=5pt, line color right=red,
inner line width right=0.5pt}};

\node[name=L, ell, shape=ellipse] {$L$};

\node[name=Yt, below right=15mm of L, ell, shape=ellipse]{$\mathbbm{1}(T_a \geq t)$};

\node[name=Ytd, right=15mm of Yt, ell, shape=ellipse]{$\mathbbm{1}(T_a \geq t + dt)$};

\node[name=A, below left =15mm of Yt, ell, shape=swig vsplit]{
\nodepart{left}{$A$}
\nodepart{right}{$a$} };

\draw[->,line width=1.0pt,>=stealth](L) to (Yt);
\draw[->,line width=1.0pt,>=stealth](L) to (Ytd);
\draw[->,line width=1.0pt,>=stealth](Yt) to (Ytd);
\draw[->,line width=1.0pt,>=stealth](A) to (Yt);
\draw[->,line width=1.0pt,>=stealth](A) to (Ytd);

\end{tikzpicture}